\documentclass[reqno]{amsart}

\usepackage{amsfonts, amsmath, amsthm, amssymb, latexsym, graphicx, geometry, xcolor, colortbl,  mathtools}

\usepackage[english]{babel}

\usepackage{physics}

\usepackage{enumitem}

\usepackage[style=numeric,doi=false,url=false,isbn=false,firstinits=true,maxbibnames=5]{biblatex}

\theoremstyle{plain}
\newtheorem{theorem}{Theorem}[section]
\newtheorem{corollary}[theorem]{Corollary}
\newtheorem{lemma}[theorem]{Lemma}

\theoremstyle{definition}

\newtheorem{remark}[theorem]{Remark}

\numberwithin{equation}{section}

\title[Free positive multiplicative Brownian motion]{Free positive multiplicative Brownian motion and the free additive convolution of semicircle and uniform distribution}

\author{Martin Auer} 
\address{Fakult\"at Mathematik, Technische Universit\"at Dortmund,
	Vogelpothsweg 87,
	D-44221 Dortmund, Germany}
\email{martin.auer@math.tu-dortmund.de}


\begin{document}
	
\date{\today}
	
\begin{abstract}
	The free positive multiplicative Brownian motion $(h_t)_{t\geq0}$ is the large $N$ limit in non-commutative distribution of matrix geometric Brownian motion.
	It can be constructed by setting $h_t:=g_{t/2}g_{t/2}^*$, where $(g_t)_{t\geq0}$ is a free multiplicative Brownian motion, which is the large $N$ limit in non-commutative distribution of the Brownian motion in $\operatorname{Gl}(N,\mathbb{C})$.
	One key property of $(h_t)_{t\geq0}$ is the fact that the corresponding spectral distributions $(\nu_t)_{t\geq0}\subset M^1((0,\infty))$ form a semigroup w.r.t.\ free multiplicative convolution.
	In recent work by M. Voit and the present author, it was shown that $\nu_t$ can be expressed by the image measure of a free {\it additive} convolution of the semicircle and the uniform distribution on an interval under the exponential map.
	In this paper, we provide a new proof of this result by calculating the moments of the free additive convolution of semicircle and uniform distributions on intervals.
	As a by-product, we also obtain new integral formulas for $\nu_t$ which generalize the corresponding known moment formulas involving Laguerre polynomials.
\end{abstract}
	
\maketitle
	
\section{Introduction}
On some filtered $W^*$-probability space $(\mathcal{A},(\mathcal{A}_t)_{t\geq0},\tau)$ let $(c_t)_{t\geq0}$ be a circular Brownian motion, see e.g.\ \cite{Nikitopoulos2022,Nikitopoulos2023} for this notion and an introduction to free stochastic calculus.
$(c_t)_{t\geq0}$ is the large $N$ limit in non-commutative distribution of the Browian motion in the space of $N\times N$ complex matrices; cf.\ \cite{Voiculescu1991} 2.3 Theorem.
%, i.e.\ $((c_t+c_t^*)/2,(c_t-c_t^*)/(2i))_{t\geq0}$ is a $2$-dimensional semicircular Brownian motion.
The {\it free $\textit{(left-)}$ multiplicative Brownian motion} $(g_t)_{t\geq0}$ is defined as the unique solution to the free stochastic differential equation (fSDE)
$$
	dg_t=g_t\,dc_t
$$
with starting condition $g_0=\mathbf{1}$ where $\mathbf{1}\in\mathcal{A}$ is the unit element.
By \cite{Kemp2016}, $(g_t)_{t\geq0}$ is the large $N$ limit in non-commutative distribution of the Brownian motion in $\operatorname{Gl}(N,\mathbb{C})$.
For positive $t$, the bounded operator $g_t$ is non-normal and its Brown measure has been fully computed in \cite{Driver2022}.
Let $(h_t:=g_{t/2}g_{t/2}^*)_{t\geq0}$ be the radial process which is called {\it free positive multiplicative Brownian motion}.
Using free stochastic calculus one can readily show that there exists a semicircular Brownian motion $(x_t)_{t\geq0}$ such that $(h_t)_{t\geq0}$ satisfies the fSDE
$$
	dh_t=\sqrt{h_t}\,dx_t\,\sqrt{h_t}+\frac{1}{2}h_t\,dt\,.
$$
Hence, $(h_t)_{t\geq0}$ can be regarded as the free probabilty analogue of geometric Brownian motion, and the distribution $\nu_t$ of $h_t$ as the free probability analogue of the $\log$-normal distribution.
The probability measure $\nu_t$ has been fully calculated in \cite{Biane1997} Section 4.2.
There, it is shown that the support of $\nu_{2t}$ is given by
$$
	\left[\left((2t+1)-2\sqrt{t(1+t)}\right)e^{-\sqrt{t(1+t)}},\left((2t+1)+2\sqrt{t(1+t)}\right)e^{\sqrt{t(1+t)}}\right]\,,
$$
and that $\nu_t$ is absolutely continuous w.r.t.\ Lebesgue measure.
The density in the interior of the support is positive and implicitly given by the solution of an inverse function problem.
The main interest in the study of $(\nu_t)_{t\geq0}$ arises from the fact that it forms a semigroup w.r.t.\ free multiplicative convolution, i.e.\ it holds $\nu_s\boxtimes\nu_t=\nu_{s+t}$ for all $s,t\geq0$.
In turn, this motivated the results in \cite{Zhong2015,Hasebe2022} concerning the support, $\log$-unimodality and density of measures of the form $\mu\boxtimes\nu_t$ for $\mu\in M^1([0,\infty))$.
Recently, it has been shown by M. Voit and the present author that the measure $\nu_t$ itself can actually be described by a free {\it additive} convolution after the change of measure under $x\mapsto\log(x)$:
\begin{theorem}[\cite{Auer2025}]\label{thm_main}
	For $t>0$ it holds
	\begin{equation}\label{eq_main_formula}
		\nu_t=\exp(\mu_{\operatorname{sc},2\sqrt{t}}\boxplus\operatorname{Unif}_{[-t/2,t/2]})\,,
	\end{equation}
	where $\exp(\cdot)$ denotes the image measure under the exponential map $x\mapsto e^x$, $\mu_{\operatorname{sc},R}$ is the semicircle distribution with radius $R>0$, and $\operatorname{Unif}_I$ is the uniform distribution on the interval $I$.
\end{theorem}
Let us briefly summarize the idea of proof:
Consider the random matrix process
$$
	(B_t+t\operatorname{diag}(\rho))_{t\geq0}\,,\;\;\rho=(-N+1,-N+3,\dots,N-3,N-1)\,,
$$
where $(B_t)_{t\geq0}$ is an $N\times N$ hermitian Brownian motion.
The authors show that the corresponding eigenvalue process is, up to a simple rescaling, equal in law to the component wise logarithm of the eigenvalue process of the matrix geometric Brownian motion $(H_t)_{t\geq0}$.
The latter matrix process is given by $H_t=G_{t/2}G_{t/2}^*$, where $(G_t)_{t\geq0}$ is the Brownian motion in $\operatorname{Gl}(N,\mathbb{C})$.
Now, on the one hand, the empirical spectral measure of $B_{t/N}+\frac{t}{N}\operatorname{diag}(\rho)$ converges to $\mu_{\operatorname{sc},2\sqrt{t}}\boxplus\operatorname{Unif}_{[-t,t]}$ for $N\to\infty$.
On the other hand, the empirical spectral measure of $H_{t/N}$ converges to $\nu_t$ for $N\to\infty$.
The result follows by combining these observations.\\
The aim of this paper is to give a new proof of Theorem \ref{thm_main} by showing equality of moments of the measures of both sides of \eqref{eq_main_formula}.
To do so, we will calculate the moments
\begin{equation}\label{eq_def_m_n}
	m_n(t):=\int_{\mathbb{R}}x^n\,d\left(\mu_{\operatorname{sc},2\sqrt{t}}\boxplus\operatorname{Unif}_{[-t,0]}\right)(x)\,.
\end{equation}
To concisely state this moment formula we use the {\it (signed) Stirling numbers of the first kind} $s(n,k)$ which are defined by the equation
$$
	x^{(n)}=\sum_{k=0}^{n}(-1)^{n-k}s(n,k)x^k\,,
$$
where $x^{(n)}=\prod_{i=0}^{n-1}(x+i)$ denotes the {\it rising factorial}.
For a broader introduction to Stirling numbers including their combinatorial interpretation, see e.g.\ \cite{Graham1989} Section 6.1.
The main part of this paper concerns the proof of the following result:
\begin{theorem}\label{thm_moment_formula}
	For $n\in\mathbb{N}_0$ and $t>0$ it holds
	$$
		m_n(t)=n!\sum_{j=\lceil n/2\rceil}^n\frac{t^j}{j!(1+j)!}s(1+j,n+1-j)\,.
	$$
\end{theorem}
Using this formula one can calculate the moments of arbitrary free additive convolutions of semicircle and uniform distributions $\mu_{\operatorname{sc},a}\boxplus\operatorname{Unif}_{[b,c]}$, see Remark \ref{rem_moments_unif_semic}.
On the other hand, by \cite{Biane1997a}, the moments of $\nu_t$ are already known and given by
%$\int_{(0,\infty)}x^n\,d\nu_t=\frac{1}{n}e^{nt/2}\sum_{j=0}^n\binom{n}{1+j}\frac{(nt)^j}{j!}$,
\begin{equation}\label{eq_moments_pos_BM_shift}
	\int_{(0,\infty)}x^n\,d\nu_t(x)
	=e^{nt/2}\frac{1}{n}L_{n-1}^{(1)}(-nt)=e^{nt/2}\frac{1}{n}\sum_{j=0}^{n-1}\binom{n}{1+j}\frac{(nt)^j}{j!}\,.
\end{equation}
Here, $L_{n}^{(\alpha)}$ denotes the $n$-th Laguerre polynomial with parameter $\alpha$, see \cite{Szegoe1975} Chapter V.
Now, Theorem \ref{thm_main} follows from Theorem \ref{thm_moment_formula} and \eqref{eq_moments_pos_BM_shift} by a straightforward computation, see Section \ref{sec_2}.
Moreover, combining Theorem \ref{thm_main} and Theorem \ref{thm_moment_formula} yields a generalization of \eqref{eq_moments_pos_BM_shift}.
To state it, recall that {\it Kummer's confluent hypergeometric function} is defined by $${}_1F_1(a;b;x):=\sum_{j=0}^{\infty}\frac{a^{(j)}}{b^{(j)}}\frac{x^j}{j!}\,.$$
We obtain the following integral formula for $\nu_t$:

\begin{corollary}\label{cor_main}
	For $t>0$ and $\alpha\in\mathbb{C}\setminus\{0\}$ it holds
	$$
		\int_{(0,\infty)}x^{\alpha}\,d\nu_t(x)
		=e^{\alpha t/2}{}_1F_1(1-\alpha;2;-\alpha t)
		=e^{\alpha t/2}\frac{1}{\alpha}\sum_{j=0}^{\infty}\binom{\alpha}{1+j}\frac{(\alpha t)^j}{j!}\,.
	$$
\end{corollary}
We will proof Theorem \ref{thm_main}, Theorem \ref{thm_moment_formula} and Corollary \ref{cor_main} in the next section.
The main work lies in proving a new combinatorial formula for the signed Stirling numbers of the first kind, see Lemma \ref{lem_stirling_identiy}.\\
We point out that the results presented in this paper are part of the author's forthcoming PhD thesis.

\section{Proof of the results}\label{sec_2}

Throughout this section, let $\nu_t$ be the distribution of the free positive multiplicative Brownian motion $h_t$ as in the introduction, and let $m_n$ be as in \eqref{eq_def_m_n}.
We first show that Theorem \ref{thm_moment_formula} and equation \eqref{eq_moments_pos_BM_shift} imply Theorem \ref{thm_main} and Corollary \ref{cor_main}:

\begin{proof}[Proof of Theorem \ref{thm_main} \& Corollary \ref{cor_main}]
	Let $\alpha\in\mathbb{C}\setminus\{0\}$.
	Then, by Theorem \ref{thm_moment_formula}, we have
	\begin{align*}
		&\int_{\mathbb{R}}e^{\alpha x}\,d\left(\mu_{\operatorname{sc},2\sqrt{t}}\boxplus\operatorname{Unif}_{[-t,0]}\right)(x)
		=\sum_{n=0}^{\infty}\frac{\alpha^n}{n!}m_n(t)\\
		=&\sum_{n=0}^{\infty}\alpha^n\sum_{j=\lceil n/2\rceil}^n\frac{1}{(1+j)!}\frac{t^j}{j!}s(1+j,n+1-j)\\
		=&-\frac{1}{\alpha}\sum_{j=0}^{\infty}\frac{(-\alpha t)^j}{(1+j)!j!}\sum_{k=1}^{1+j}(-1)^{j+1-k}s(1+j,k)(-\alpha)^k\\
		=&-\frac{1}{\alpha}\sum_{j=0}^{\infty}\frac{(-\alpha)^{(j+1)}}{(1+j)!}\frac{(-\alpha t)^j}{j!}
		={}_1F_1(1-\alpha;2;-\alpha t)\,.
	\end{align*}	
	In particular, it holds for $n\in\mathbb{N}$:
	\begin{align*}
		\int_{\mathbb{R}}x^n\,d\exp\left(\mu_{\operatorname{sc},2\sqrt{t}}\boxplus\operatorname{Unif}_{[-t,0]}\right)(x)
		={}_1F_1(1-n;2;-nt)
		=\frac{1}{n}\sum_{j=0}^{n-1}\binom{n}{1+j}\frac{(nt)^j}{j!}\,.
	\end{align*}
	By \eqref{eq_moments_pos_BM_shift}, this is the $n$-th moment of $D_{e^{-t/2}}(\nu_t)$.
	Since measures with compact support are uniquely characterized by their moments, this shows
	\begin{equation*}
		\exp\left(\mu_{\operatorname{sc,2\sqrt{t}}}\boxplus\operatorname{Unif}_{[-t,0]}\right)
		=D_{e^{-t/2}}(\nu_t)\,.
	\end{equation*}
	Now, scaling in space by the factor $e^{t/2}$ gives Theorem \ref{thm_main}.
	In turn also Corollary \ref{cor_main} holds:
	\begin{align*}
		\int_{(0,\infty)}x^{\alpha}\,d\nu_t
		=\int_{\mathbb{R}}e^{\alpha x}\,d\left(\mu_{\operatorname{sc},2\sqrt{t}}\boxplus\operatorname{Unif}_{[-t/2,t/2]}\right)(x)
		=&e^{\alpha t/2}\int_{\mathbb{R}}e^{\alpha x}\,d\left(\mu_{\operatorname{sc},2\sqrt{t}}\boxplus\operatorname{Unif}_{[-t,0]}\right)(x)\\
		=&e^{\alpha t/2}{}_1F_1(1-\alpha;2;-\alpha t)\,.
	\end{align*}
\end{proof}

Hence, it remains to prove Theorem \ref{thm_moment_formula}.
Before we state the proof, we note another straightforward consequence of this result:

\begin{remark}\label{rem_moments_unif_semic}
	Using Theorem \ref{thm_moment_formula} we can calculate the moments of general free additive convolutions of semicircle and uniform distributions on intervals as follows: Let $a>0$ and $b,c\in\mathbb{R}$ with $b<c$. Then
	$$
	\mu_{\operatorname{sc},2\sqrt{a}}\boxplus\operatorname{Unif}_{[b,c]}
	=\mu_{\operatorname{sc},2\sqrt{a}}\boxplus\operatorname{Unif}_{[b-c,0]}\boxplus\delta_{c}\,,
	$$
	where we used that classical additive and free additive convolutions with Dirac measures coincide, and that the free additive convolution is associative.
	For $\gamma\in\mathbb{R}$ let
	\begin{equation}\label{eq_scaling_map}
		D_{\gamma}\colon\mathbb{R}\to\mathbb{R},\,x\mapsto\gamma x
	\end{equation}
	be the scaling map with factor $\gamma$.
	Then, the image measure of a free additive convolution under $D_{\gamma}$ linearizes as in the classical case: $D_{\gamma}(\mu\boxplus\nu)=D_{\gamma}(\mu)\boxplus D_{\gamma}(\nu)$, $\mu,\nu\in M^1(\mathbb{R})$.
	Thus
	\begin{align*}
		\mu_{\operatorname{sc},2\sqrt{a}}\boxplus\operatorname{Unif}_{[b,c]}
		=&\mu_{\operatorname{sc},2\sqrt{a}}\boxplus\operatorname{Unif}_{[b-c,0]}\boxplus\delta_{c}\\
		=&D_{a/(c-b)}(\mu_{\operatorname{sc},2\sqrt{(c-b)^2/a}})\boxplus D_{a/(c-b)}(\operatorname{Unif}_{[-(c-b)^2/a,0]})\boxplus\delta_{c}\\
		=&D_{a/(c-b)}\left(\left.\left(\mu_{\operatorname{sc},2\sqrt{t}})\boxplus \operatorname{Unif}_{[-t,0]}\right)\right\rvert_{t=(c-b)^2/a}\right)\ast\delta_{c}\,.
	\end{align*}
	Here $\ast$ denotes the classical (additive) convolution of real measures.
	Using the binomial theorem we arrive at
	$$
	\int_{\mathbb{R}}x^n\,d\left(\mu_{\operatorname{sc},2\sqrt{a}}\boxplus\operatorname{Unif}_{[b,c]}\right)(x)
	=n!\sum_{k=0}^{n}\frac{c^{n-k}}{(n-k)!}\sum_{j=\lceil k/2\rceil}^k\frac{(c-b)^{2j-k}a^{k-j}}{j!(1+j)!}s(1+j,k+1-j)\,.
	$$
	In particular, this explains why we chose to define $(m_n(t))_{n\geq0}$ as the moments of $\mu_{\operatorname{sc},2\sqrt{t}}\boxplus\operatorname{Unif}_{[-t,0]}$ and not as the moments of $\mu_{\operatorname{sc},2\sqrt{t}}\boxplus\operatorname{Unif}_{[-t/2,t/2]}$.
\end{remark}

In order to show Theorem \ref{thm_moment_formula}, we will apply the following two Lemmas.
In the corresponding proofs, the following notation is used:
Let $f(z)=\sum_{n=0}^{\infty}a_nz^n$ be a power series. We will denote the $n$-th coefficient by $[z^n]f(z):=a_n$ and call $[z^n]$ {\it coefficient extractor}.\\
The first of the following two Lemmas states a time inhomogeneous recursion ODE for the moments $m_n$:

\begin{lemma}\label{lem_unif_semic_moment_recur}
	$m_n$ is a polynomial of degree $n$ with leading coefficient $[t^n]m_n(t)=(-1)^n/(1+n)$, $n\in\mathbb{N}_0$.
	Further, it holds for all $t>0$:
	$$
	\frac{d}{dt}m_n(t)=n\left(t^{-1}m_n(t)-\frac{1}{2}\sum_{j=0}^{n-2}m_j(t)m_{n-j-2}(t)\right)\,.
	$$
\end{lemma}

\begin{proof}
	On some $W^*$-probability space $(\mathcal{A},\tau)$, let $(x_t)_{t\geq0}$ be a semicircular Brownian motion and $a$ a self-adjoint bounded operator which is freely independent of $(x_t)_{t\geq0}$ and whose spectral distribution equals $\operatorname{Unif}_{[-1,0]}$.
	Denote by $\tilde{m}_n(t)$ the $n$-th moment of
	$$
	\mu_t:=\operatorname{Unif}_{[-1,0]}\boxplus\mu_{\operatorname{sc},2\sqrt{t}}\,.
	$$
	By construction we have $\tilde{m}_n(t)=\tau((a+x_t)^n)$.
	Hence, by \cite{Nikitopoulos2022} Theorem 3.5.3, it holds
	$$
	\frac{d}{dt}\tilde{m}_n(t)=\frac{n}{2}\sum_{j=0}^{n-2}\tilde{m}_j(t)\tilde{m}_{n-j-2}(t)\,.
	$$
	Induction on $n$ shows that $\tilde{m}_n$ is a polynomial of degree $n-1$.
	Note that it holds $[t^0]\tilde{m}_n=\tilde{m}_n(0)=(-1)^n/(1+n)$.
	Let $D_{\gamma}$ be as in \eqref{eq_scaling_map}.
	Then, using $D_{t}(\mu_{1/t})=\operatorname{Unif}_{[-t,0]}\boxplus\mu_{\operatorname{sc},2\sqrt{t}}$, we have $m_n(t)=t^n\tilde{m}_n(1/t)$ for $t>0$.
	This shows that $m_n$ is a polynomial of degree $n$ with $[t^n]m_n(t)=(-1)^n/(1+n)$.
	Moreover, we have
	\begin{align*}
		\frac{d}{dt}m_n(t)
		=&nt^{n-1}\tilde{m}_n(1/t)-\frac{n}{2}t^{n-2}\sum_{j=0}^{n-2}\tilde{m}_j(1/t)\tilde{m}_{n-2-j}(1/t)\\
		=&nt^{-1}m_n(t)-\frac{n}{2}\sum_{j=0}^{n-2}m_j(t)m_{n-2-j}(t)\,.
	\end{align*}
\end{proof}

The second Lemma is a combinatorial identity involving a double sum over the product of two Stirling numbers and three binomial coefficients, one of which is inverted.
We postpone its technical proof to the end of this section. 

\begin{lemma}\label{lem_stirling_identiy}
	For all $l,m\in\mathbb{N}$ it holds
	\begin{equation}\label{stirling_number_double_sum_eq}
		\begin{split}
			&2ls(1+m,1+l)\\
			=&\frac{m+1}{l+m-1}\sum_{n=1}^l\sum_{k=0}^{m-1}\binom{l+m-2}{n+k-1}^{-1}\binom{m}{k}\binom{m}{1+k}s(1+k,n)s(m-k,l+1-n)\,.
		\end{split}
	\end{equation}
\end{lemma}

\begin{proof}[Proof of Theorem \ref{thm_moment_formula}]
	Set $c(n,k):=[t^k]m_n(t)$, $n\in\mathbb{N}$, $k\in\{0,\dots,n\}$.
	Then, by Lemma \ref{lem_unif_semic_moment_recur}, we know that for all $n\geq1$ we have $c(n,0)=0$, $c(n,n)=(-1)^n/(1+n)=n!\frac{1}{n!(1+n)!}s(1+n,n+1-n)$ and:
	\begin{gather*}
		\sum_{k=1}^{n}kc(n,k)t^{k-1}
		=n\sum_{k=1}^nc(n,k)t^{k-1}-\frac{n}{2}\sum_{j=0}^{n-2}\sum_{l=0}^{j}c(j,l)t^l\sum_{m=0}^{n-j-2}c(n-j-2,m)t^m\,.
	\end{gather*}
	We can rearrange this in the following way:
	\begin{align*}
		\sum_{k=1}^n(n-k)c(n,k)t^{k-1}
		=&\frac{n}{2}\sum_{l=0}^{n-2}\sum_{j=l}^{n-2}\sum_{m=0}^{n-j-2}t^{m+l}c(j,l)c(n-j-2,m)\\
		=&\frac{n}{2}\sum_{k=1}^{n-1}t^{k-1}\sum_{l=0}^{k-1}\sum_{j=l}^{l+n-1-k}c(j,l)c(n-j-2,k-1-l)\,,
	\end{align*}
	where we set $k:=l+m+1$ in the second equality.
	Comparison of coefficients yields
	\begin{align*}
		(n-k)c(n,k)=&\frac{n}{2}\sum_{l=0}^{k-1}\sum_{j=l}^{l+n-1-k}c(j,l)c(n-j-2,k-1-l)\\
		=&\frac{n}{2}\sum_{l=0}^{k-1}\sum_{j=1}^{n-k}c(j+l-1,l)c(n-l-j-1,k-1-l)\,.
	\end{align*}
	Now, we show $c(n,k)=n!\frac{1}{k!(1+k)!}s(1+k,n+1-k)$, $k\in\{0,\dots,n\}$, by induction on $n\in\mathbb{N}_0$.
	We already know that this formula holds true for $c(n,n)$ and $c(n,0)$, $n\in\mathbb{N}$.
	The induction start at $n=0$ holds since $c(0,0)=1=s(1,1)$ and the induction start at $n=1$ holds by the already proven formulas for $c(n,n)$ and $c(n,0)$.
	Suppose $n\geq2$ and that the induction assertion has been proven for all $n'\in\{0,\dots,n-2\}$.
	Then for $k\in\{1,\dots,n-1\}$, by induction assumption and Lemma \ref{lem_stirling_identiy}, we have:
	\begin{align*}
		&c(n,k)\\
		=&\frac{1}{2}\frac{n}{n-k}\sum_{l=0}^{k-1}\sum_{j=1}^{n-k}c(j+l-1,l)c(n-l-j-1,k-1-l)\\
		=&\frac{1}{2}\frac{n}{n-k}\sum_{l=0}^{k-1}\sum_{j=1}^{n-k}\frac{(j+l-1)!}{l!(1+l)!}s(1+l,j)\frac{(n-l-j-1)!}{(k-1-l)!(k-l)!}s(k-l,n+1-j-k)\\
		=&\frac{1}{2}\frac{n}{n-k}\frac{1}{k!k!}(n-2)!\sum_{l=0}^{k-1}\sum_{j=1}^{n-k}\binom{n-2}{j+l-1}^{-1}\binom{k}{l}\binom{k}{1+l}s(1+l,j)s(k-l,n+1-j-k)\\
		=&\frac{1}{2}\frac{n}{n-k}\frac{1}{k!k!}\frac{(n-1)!}{1+k}2(n-k)s(1+k,1+n-k)=n!\frac{1}{k!(1+k)!}s(1+k,n+1-k)\,.
	\end{align*}
	Lastly, it holds $c(n,k)=0$ for all $k<\lceil n/2\rceil$ since we have $1+k<n+1-k$, which implies $s(1+k,n+1-k)=0$.
\end{proof}

In the remainder of this paper we prove Lemma \ref{lem_stirling_identiy}, which is the last missing link in the proof of Theorem \ref{thm_moment_formula}.
We will use Egorychev's method which provides a variety of effective tools to show combinatorial formulas (cf. \cite{Egorychev1984,Riedel2023}).

\begin{proof}[Proof of Lemma \ref{lem_stirling_identiy}]
	We start by showing \eqref{stirling_number_double_sum_eq} in the easiest case $m=1$:
	$$
	2ls(2,1+l)
	=\begin{rcases}
		\begin{dcases}
			2\,,&l=1\\
			0\,,&\text{else}
		\end{dcases}
	\end{rcases}
	=\frac{2}{l}s(1,l)=\frac{2}{l}\sum_{n=1}^l\binom{l-1}{n-1}^{-1}s(1,n)s(1,l+1-n)\,.
	$$
	Now let $m\geq2$. Rearranging the double sum of the right hand side of \eqref{stirling_number_double_sum_eq} gives
	\begin{equation}\label{pf_stirling_recur_eq0}
		\begin{split}
			&\frac{m+1}{l+m-1}\sum_{n=1}^l\sum_{k=0}^{m-1}\binom{l+m-2}{n+k-1}^{-1}\binom{m}{k}\binom{m}{1+k}s(1+k,n)s(m-k,l+1-n)\\
			=&\frac{m+1}{l+m-1}\sum_{n=0}^{l+1}\sum_{k=1}^{m-2}\binom{l+m-2}{n+k-1}^{-1}\binom{m}{k}\binom{m}{1+k}s(1+k,n)s(m-k,l+1-n)\\
			&+2\frac{m(m+1)}{l+m-1}s(m,l)\\
			=&:A+2\frac{m(m+1)}{l+m-1}s(m,l)\,.
		\end{split}
	\end{equation}
	Now, we apply techniques of the Egorychev method, i.e.\ we transform a suitable generating series of $A$.
	To this end we make use of the identities
	\begin{gather*}
		s(n,k)=\frac{n!}{k!}\left[z^n\right]\log^k(1+z)\,,\quad n,k\in\mathbb{N}_0\,,\\
		\int_0^1t^n(1-t)^k\,dt=\frac{1}{n+k-1}\binom{n+k}{n}^{-1}\,,\quad n,k\in\mathbb{N}_0\,,
	\end{gather*}
	see \cite{Riedel2023} Section 5 for the first identity, the second identity is a special case of the well known beta integral.
	Now, we can write $A$ in the following way:
	\begin{align*}
		A
		=&(m+1)\sum_{k=1}^{m-2}\sum_{n=0}^{l+1}\int_0^1t^{n+k-1}(1-t)^{l+m-1-n-k}\,dt\,\binom{m}{k}\binom{m}{1+k}\\
		&\cdot\left[x^{1+k}\right]\left[y^{m-k}\right]\frac{(1+k)!}{n!}\frac{(m-k)!}{(l+1-n)!}\log^n(1+x)\log^{l+1-n}(1+y)\\
		=&\frac{m(m+1)!}{(l+1)!}\sum_{k=1}^{m-2}\binom{m-1}{k}\left[x^{1+k}\right]\left[y^{m-k}\right]\\
		&\cdot\int_0^1t^{k-1}(1-t)^{m-2-k}\sum_{n=0}^{l+1}\binom{l+1}{n}\left(t\log(1+x)\right)^n\left((1-t)\log(1+y)\right)^{l+1-n}\,dt\\
		=&\frac{m(m+1)!}{(l+1)!}\sum_{k=1}^{m-2}\binom{m-1}{k}\left[x^{1+k}\right]\left[y^{m-k}\right]\\
		&\cdot\int_0^1t^{k-1}(1-t)^{m-2-k}\left(t\log(1+x)+(1-t)\log(1+y)\right)^{l+1}\,dt\,.
	\end{align*}
	Using coefficient extraction for the exponential series, we can rewrite this further as
	\begin{align*}
		A=&m(m+1)!\left[\alpha^{l+1}\right]\sum_{k=1}^{m-2}\binom{m-1}{k}\left[x^{1+k}\right]\left[y^{m-k}\right]\\
		&\cdot\int_0^1t^{k-1}(1-t)^{m-2-k}(1+y)^{\alpha}\exp\left(\log\left(\frac{1+x}{1+y}\right)\alpha t\right)\,dt\,.
	\end{align*}
	Now, we can apply Euler's integral transform for Kummer's confluent hypergeometric function ${}_1F_1$ which states that
	\begin{equation*}
		{}_1F_1(a;b;x)=\frac{\Gamma(b)}{\Gamma(a)\Gamma(b-a)}\int_0^1t^{a-1}(1-t)^{b-a-1}e^{tx}\,dt\,.
	\end{equation*}
	This gives
	\begin{align*}
		&\frac{A}{m(m+1)!}\\
		=&\left[\alpha^{l+1}\right]\sum_{k=1}^{m-2}\binom{m-1}{k}\left[x^{1+k}\right]\left[y^{m-k}\right](1+y)^{\alpha}\frac{(k-1)!(m-2-k)!}{(m-2)!}\\
		&\cdot{}_1F_1\left(k;m-1;\log\left(\frac{1+x}{1+y}\right)\alpha\right)\\
		=&\left[\alpha^{l+1}\right]\sum_{k=1}^{m-2}\frac{m-1}{k(m-1-k)}\left[x^{1+k}\right]\left[y^{m-k}\right](1+y)^{\alpha}{}_1F_1\left(k;m-1;\alpha\log\left(\frac{1+x}{1+y}\right)\right)\\
		=&\left[\alpha^{l+1}\right]\sum_{k=1}^{m-2}\left(\frac{1}{k}+\frac{1}{m-1-k}\right)\left[x^{1+k}\right]\left[y^{m-k}\right](1+y)^{\alpha}{}_1F_1\left(k;m-1;\alpha\log\left(\frac{1+x}{1+y}\right)\right)\,.
	\end{align*}
	Reversing the order of summation in one part of the sum and applying Kummer's transformation rule ${}_1F_1(a;b;x)=e^x{}_1F_1(b-a;b;-x)$, we see
	\begin{align*}
		&\sum_{k=1}^{m-2}\frac{1}{m-1-k}\left[x^{1+k}\right]\left[y^{m-k}\right](1+y)^{\alpha}{}_1F_1\left(k;m-1;\alpha\log\left(\frac{1+x}{1+y}\right)\right)\\
		=&\sum_{k=1}^{m-2}\frac{1}{k}\left[x^{m-k}\right]\left[y^{1+k}\right](1+y)^{\alpha}{}_1F_1\left(m-1-k;m-1;\alpha\log\left(\frac{1+x}{1+y}\right)\right)\\
		=&\sum_{k=1}^{m-2}\frac{1}{k}\left[x^{m-k}\right]\left[y^{1+k}\right](1+x)^{\alpha}{}_1F_1\left(k;m-1;\alpha\log\left(\frac{1+y}{1+x}\right)\right)\,.
	\end{align*}
	Thus
	\begin{align*}
		\frac{A}{m(m+1)!}
		=&2\left[\alpha^{l+1}\right]\sum_{k=1}^{m-2}\frac{1}{k}\left[x^{1+k}\right]\left[y^{m-k}\right](1+y)^{\alpha}{}_1F_1\left(k;m-1;\alpha\log\left(\frac{1+x}{1+y}\right)\right)\,.
	\end{align*}
	Now, we can extend the sum over $k$ to infinity with the following correction due to $[y^{m-k}]$:
	\begin{align*}
		\frac{A}{2m(m+1)!}
		=&\left[\alpha^{l+1}\right]\sum_{k=1}^{\infty}\frac{1}{k}\left[x^{1+k}\right]\left[y^{m-k}\right](1+y)^{\alpha}{}_1F_1\left(k;m-1;\alpha\log\left(\frac{1+x}{1+y}\right)\right)\\
		&-\left[\alpha^{l+1}\right]\frac{1}{m-1}\left[x^{m}\right]\left[y^{1}\right](1+y)^{\alpha}{}_1F_1\left(m-1;m-1;\alpha\log\left(\frac{1+x}{1+y}\right)\right)\\
		&-\left[\alpha^{l+1}\right]\frac{1}{m}\left[x^{1+m}\right]\left[y^{0}\right](1+y)^{\alpha}{}_1F_1\left(m;m-1;\alpha\log\left(\frac{1+x}{1+y}\right)\right)\\
		=&:[\alpha^{l+1}](S_1+S_2+S_3)\,.
	\end{align*}
	We will calculate the three appearing terms seperately, where the main work lies in the calculation of $S_1$.
	For the latter task we aim to identify a geometric series in the summation over $k$.
	For this, we use Cauchy's integral formula in the following way:
	\begin{align*}
		S_1
		=&\sum_{k=1}^{\infty}\frac{1}{k}\left[x^{1+k}\right]\left[y^{m-k}\right](1+y)^{\alpha}\sum_{j=0}^{\infty}\frac{k^{(j)}\alpha^j}{(m-1)^{(j)}j!}\log^j\left(\frac{1+x}{1+y}\right)\\
		=&\sum_{k=1}^{\infty}\frac{1}{k}\left(\frac{1}{2\pi i}\right)^2\int_{\lvert x\rvert=1/2}\int_{\lvert y\rvert=1/4}x^{-(k+2)}y^{k-(m+1)}(1+y)^{\alpha}\\
		&\cdot\sum_{j=0}^{\infty}\frac{k^{(j)}\alpha^j}{(m-1)^{(j)}j!}\log^j\left(\frac{1+x}{1+y}\right)\,dydx\\
		=&\left(\frac{1}{2\pi i}\right)^2\int_{\lvert y\rvert=1/4}\int_{\lvert x\rvert=1/2}x^{-2}y^{-(m+1)}(1+y)^{\alpha}\sum_{j=1}^{\infty}\frac{\alpha^j}{(m-1)^{(j)}j!}\log^j\left(\frac{1+x}{1+y}\right)\\
		&\cdot\sum_{k=1}^{\infty}\left.\frac{d^{j-1}}{dz^{j-1}}z^{k+j-1}\right\rvert_{z=\frac{y}{x}}\,dxdy\,.
	\end{align*}
	Now, one can easily verify the following identity by induction over $j$: 
	\begin{align*}
		\sum_{k=1}^{\infty}\frac{d^{j-1}}{dz^{j-1}}z^{k+j-1}=\frac{d^{j-1}}{dz^{j-1}}\frac{z^j}{1-z}=(j-1)!\left((1-z)^{-j}-1\right)\,,\quad \lvert z\rvert<1\,.
	\end{align*}
	Plugging this in we receive
	\begin{align*}
		S_1
		=&\left(\frac{1}{2\pi i}\right)^2\int_{\lvert y\rvert=1/4}\int_{\lvert x\rvert=1/2}x^{-2}y^{-(m+1)}(1+y)^{\alpha}\\
		&\cdot\sum_{j=1}^{\infty}\frac{\alpha^j}{(m-1)^{(j)}j}\log^j\left(\frac{1+x}{1+y}\right)\left(\left(\frac{x}{x-y}\right)^j-1\right)\,dxdy\,.
	\end{align*}
	Note that, by using $\log(1+z)=-\sum_{k=1}^{\infty}\frac{(-z)^k}{k}$, $\lvert z\rvert<1$, for all $j\geq2$ and $0<\lvert y\rvert<1/2$ we have
	\begin{align*}
		&\int_{\lvert x\rvert=1/2}x^{j-2}\log^j\left(\frac{1+x}{1+y}\right)\left(\frac{1}{x-y}\right)^j\,dx\\
		=&\int_{\lvert x\rvert=1/2}(x+y)^{j-2}\log^j\left(1+\frac{x}{1+y}\right)\left(\frac{1}{x}\right)^j\,dx
		=0\,.
	\end{align*}
	This gives
	\begin{align*}
		&\left(2\pi i\right)^2S_1\\
		=&\int_{\lvert y\rvert=1/4}\int_{\lvert x\rvert=1/2}x^{-1}(x-y)^{-1}y^{-(m+1)}(1+y)^{\alpha}\frac{\alpha}{m-1}\log\left(\frac{1+x}{1+y}\right)\,dxdy\\
		&-\int_{\lvert y\rvert=1/2}\int_{\lvert x\rvert=1/2}x^{-2}y^{-(m+1)}(1+y)^{\alpha}\sum_{j=1}^{\infty}\frac{\alpha^j}{(m-1)^{(j)}j}\log^j\left(\frac{1+x}{1+y}\right)\,dxdy\\
		=&:\left(2\pi i\right)^2\left(S_{1,1}+S_{1,2}\right)\,.
	\end{align*}
	To simplify $S_{1,1}$ we first observe
	\begin{align*}
		&\int_{\lvert x\rvert=1/2}x^{-1}(x-y)^{-1}\log\left(\frac{1+x}{1+y}\right)\,dx\\
		=&-y^{-1}\int_{\lvert x\rvert=1/2}x^{-1}\log\left(\frac{1+x}{1+y}\right)\,dx
		+y^{-1}\int_{\lvert x\rvert=1/2}(x-y)^{-1}\log\left(\frac{1+x}{1+y}\right)\,dx\\
		=&(2\pi i)y^{-1}\log(1+y)+y^{-1}\int_{\lvert x\rvert=1/2}x^{-1}\log\left(1+\frac{x}{1+y}\right)\,dx
		=(2\pi i)y^{-1}\log(1+y)\,.
	\end{align*}
	This gives
	\begin{equation}\label{pf_stirling_recur_eq1}
		\begin{split}
			S_{1,1}
			=&\frac{1}{2\pi i}\int_{\lvert y\rvert=1/4}y^{-(m+2)}(1+y)^{\alpha}\frac{\alpha}{m-1}\log(1+y)\,dy\\
			%	 		=&-\frac{1}{m-1}\left[x^0\right]\left[y^{m+1}\right](1+y)^{\alpha}\alpha(\log(1+x)-\log(1+y))\\
			=&\frac{1}{m-1}\left[y^{m+1}\right](1+y)^{\alpha}\log\left((1+y)^{\alpha}\right)\\
			=&\frac{1}{m-1}\frac{1}{(m+1)!}\sum_{j=0}^{\infty}(j+1)\alpha^{j+1}s(m+1,j+1)\,.
		\end{split}
	\end{equation}
	
	The second contribution to $S_1$ can be simplified in the following way:
	
	\begin{align*}
		S_{1,2}
		=&-\left(\frac{1}{2\pi i}\right)^2\int_{\lvert y\rvert=1/2}\int_{\lvert x\rvert=1/2}x^{-2}y^{-(m+1)}(1+y)^{\alpha}\sum_{j=1}^{\infty}\frac{\alpha^j}{(m-1)^{(j)}j}\log^j\left(\frac{1+x}{1+y}\right)\,dxdy\\
		=&-\left[x^1\right]\left[y^m\right](1+y)^{\alpha}\sum_{j=1}^{\infty}\frac{\alpha^j}{(m-1)^{(j)}j}\sum_{k=0}^j\binom{j}{k}\log^k(1+x)\left(-\log(1+y)\right)^{j-k}\\
		=&-\left[y^m\right](1+y)^{\alpha}\sum_{j=1}^{\infty}\frac{\alpha^j(-1)^{j-1}}{(m-1)^{(j)}}\log^{j-1}(1+y)\\
		=&-\left[y^m\right]\sum_{n=0}^{\infty}\frac{\alpha^n}{n!}\sum_{j=1}^{\infty}\frac{\alpha^j(-1)^{j-1}(m-2)!}{(m+j-2)!}\log^{n+j-1}(1+y)\\
		=&-(m-2)!\left[y^m\right]\sum_{k=0}^{\infty}\alpha^{k+1}\log^k(1+y)\sum_{n=0}^k\frac{(-1)^{n-k}}{n!(m+k-n-1)!}\\
		=&\frac{(m-2)!}{m!}\sum_{k=0}^{\infty}(-\alpha)^{k+1}k!s(m,k)\frac{1}{(m+k-1)!}\sum_{n=0}^k(-1)^n\binom{m+k-1}{n}\\
		=&\frac{(m-2)!}{m!}\sum_{k=0}^{\infty}(-\alpha)^{k+1}k!s(m,k)\frac{(-1)^k}{(m+k-1)!}\binom{m+k-2}{k}\\
		=&-\frac{1}{m!}\sum_{k=0}^{\infty}\frac{\alpha^{k+1}}{m+k-1}s(m,k)\,,
	\end{align*}
	where we used the identity $\sum_{n=0}^k(-1)^n\binom{N}{i}=(-1)^k\binom{N-1}{k}$ in the second last equality.
	This completes the calculation of $S_1$.\\
	Using ${}_1F_1(m-1;m-1;x)=e^x$, we see that $S_2$ does not contribute:
	\begin{align*}
		S_2
		=&-\frac{1}{m-1}\left[x^m\right]\left[y^1\right](1+y)^{\alpha}\left(\frac{1+x}{1+y}\right)^{\alpha}
		=0\,.
	\end{align*}
	Applying ${}_1F_1(m;m-1;x)=e^x{}_1F_1(-1;m-1;-x)=e^x(1+x/(m-1))$, we can simplify the remaining sum as follows:
	\begin{align*}
		S_3
		=&-\frac{1}{m}\left[x^{m+1}\right]\left[y^0\right](1+y)^{\alpha}\left(\frac{1+x}{1+y}\right)^{\alpha}\left(1+\frac{\alpha}{m-1}\log\left(\frac{1+x}{1+y}\right)\right)\\
		=&-\frac{1}{m}\left[x^{m+1}\right](1+x)^{\alpha}\left(1+\frac{\alpha}{m-1}\log(1+x)\right)\\
		=&-\frac{1}{m(m+1)!}\sum_{j=0}^{\infty}\alpha^js(m+1,j)-\frac{1}{m(m-1)(m+1)!}\sum_{j=0}^{\infty}(j+1)\alpha^{j+1}s(m+1,j+1)\,.
	\end{align*}
	Combining the previous calculations we get
	\begin{align*}
		A
		=&2m(m+1)!\left[\alpha^{l+1}\right](S_{1,1}+S_{1,2}+S_2+S_3)\\
		=&2m(m+1)!\left(\frac{l+1}{(m-1)(m+1)!}s(m+1,l+1)-\frac{1}{m!(m+l-1)}s(m,l)\right.\\
		&\left.-\frac{1}{m(m+1)!}s(m+1,l+1)-\frac{l+1}{m(m-1)(m+1)!}s(m+1,l+1)\right)\\
		=&-2\frac{m(m+1)}{m+l-1}s(m,l)+2m\frac{(l+1)m-(m-1)-(l+1)}{m(m-1)}s(m+1,l+1)\\
		=&-2\frac{m(m+1)}{m+l-1}s(m,l)+2ls(m+1,l+1)\,.
	\end{align*}
	Plugging this back into \eqref{pf_stirling_recur_eq0} shows the claim.
\end{proof}

\printbibliography

\end{document}